\documentclass[12pt]{amsart}
\usepackage{amssymb}

\newcommand{\Arh}{Arhangel'ski\u{\i}}

\long\def\forget#1\forgotten{}
\newcommand{\issuenumber}{19}
\newcommand{\issuemonth}{December}
\newcommand{\issueyear}{2006}

\newtheorem{thm}{Theorem}[section]
\newtheorem{prob}[thm]{Problem}

\newtheorem{issue}{Issue}

\theoremstyle{definition}

\theoremstyle{remark}

\newcommand{\ed}{
\newpage

\section{Unsolved problems from earlier issues}

\begin{issue}
Is $\binom{\Omega}{\Gamma}=\binom{\Omega}{\Tau}$?
\end{issue}

\begin{issue}
Is $\ufin(\cO,\Omega)=\sfin(\Gamma,\Omega)$?
And if not, does $\ufin(\cO,\Gamma)$ imply
$\sfin(\Gamma,\Omega)$?
\end{issue}

\stepcounter{issue}

\begin{issue}
Does $\sone(\Omega,\Tau)$ imply $\ufin(\Gamma,\Gamma)$?
\end{issue}

\begin{issue}
Is $\fp=\fp^*$? (See the definition of $\fp^*$ in that issue.)
\end{issue}

\begin{issue}
Does there exist (in ZFC) an uncountable set satisfying $\sone(\BG,\B)$?
\end{issue}

\stepcounter{issue}

\begin{issue}
Does $X \nin \NON(\M)$ and $Y\nin\mathsf{D}$ imply that
$X\cup Y\nin \COF(\M)$?
\end{issue}

\begin{issue}[CH]
Is $\split(\Lambda,\Lambda)$ preserved under finite unions?
\end{issue}

\begin{issue}
Is $\cov(\M)=\fo$? (See the definition of $\fo$ in that issue.)
\end{issue}

\begin{issue}
Does $\sone(\Gamma,\Gamma)$ always contain an element of cardinality $\fb$?
\end{issue}

\begin{issue}
Could there be a Baire metric space $M$ of weight $\aleph_1$ and a partition
$\mathcal{U}$ of $M$ into $\aleph_1$ meager sets where for each ${\mathcal U}'\subset\mathcal U$,
$\bigcup {\mathcal U}'$ has the Baire property in $M$?
\end{issue}

\stepcounter{issue} 

\begin{issue}
Does there exist (in ZFC) a set of reals $X$ of cardinality $\fd$ such that all
finite powers of $X$ have Menger's property $\sfin(\cO,\cO)$?
\end{issue}

\begin{issue}
Can a Borel non-$\sigma$-compact group be generated by a Hurewicz subspace?
\end{issue}

\begin{issue}[MA]
Is there an uncountable $X\sbst\R$ satisfying $\sone(\BO,\BG)$?
\end{issue}

\begin{issue}[CH]
Is there a totally imperfect $X$ satisfying $\ufin(\cO,\Gamma)$
that can be mapped continuously onto $\Cantor$?
\end{issue}

\begin{issue}[CH]
Is there a Hurewicz $X$ such that $X^2$ is Menger but not Hurewicz?
\end{issue}

\begin{issue}
Does the Pytkeev property of $C_p(X)$ imply the Menger property of $X$?
\end{issue}

\general\end{document}}

\newcommand{\Cantor}{{\{0,1\}^\N}}

\newcommand{\fb}{\mathfrak{b}}
\newcommand{\fg}{\mathfrak{g}}
\newcommand{\fc}{\mathfrak{c}}
\newcommand{\fd}{\mathfrak{d}}
\newcommand{\fp}{\mathfrak{p}}

\newcommand{\NON}{{\mathsf   {NON}}}
\newcommand{\COF}{{\mathsf   {COF}}}

\newcommand{\M}{\mathcal{M}}

\newcommand{\cov}{\mathsf{cov}}

\newcommand{\R}{\mathbb{R}}

\newcommand{\fo}{\mathfrak{od}}
\newcommand{\cl}[1]{\overline{#1}}

\newcommand{\w}{\omega}

\renewcommand{\split}{\mathsf{Split}}
\newcommand{\bq}{\begin{quote}}
\newcommand{\eq}{\end{quote}}
\newcommand{\cO}{\mathcal{O}}
\newcommand{\B}{\mathcal{B}}
\newcommand{\BG}{\B_\Gamma}

\newcommand{\BO}{\B_\Omega}

\newcommand{\sone}{\mathsf{S}_1}    \newcommand{\sfin}{\mathsf{S}_{fin}}

\newcommand{\ufin}{\mathsf{U}_{fin}}
    
\newcommand{\seq}[1]{\{#1\}_{n\in\N}}

\newcommand{\nin}{\not\in}

\newcommand{\cU}{\mathcal{U}}

\newcommand{\NN}{{\N^\N}}
\newcommand{\N}{\mathbb{N}}
\newcommand{\Z}{\mathbb{Z}}

\newcommand{\sbst}{\subseteq}
\newcommand{\by}[2]{\par\hfill\emph{#1}, #2}
\newcommand{\nby}[1]{\par\hfill\emph{#1}}
\newcommand{\Tau}{\mathrm{T}}
\newcommand{\CE}{\textsc{CE}}

\newcommand{\be}{\begin{enumerate}}
\newcommand{\ee}{\end{enumerate}}
\newcommand{\bi}{\begin{itemize}}
\newcommand{\ei}{\end{itemize}}
\newcommand{\itm}{\item}

\newcommand{\general}{\small\vfill\par\noindent\hrulefill\par
\noindent\textbf{Previous issues.} The previous issues of this
bulletin, which contain general information (first issue), basic
definitions, research announcements, and open problems (all
issues) are available online,
at \texttt{http://front.math.ucdavis.edu/search?\&t=\%22SPM+Bulletin\%22}
\\[0.1cm]
\textbf{Contributions.}
Please submit your contributions (announcements, discussions, and open problems)
by e-mailing us. It is preferred to write them
in \LaTeX{}.
The authors are urged to use as standard notation as possible, or otherwise give
the definitions or a reference to where the notation is explained.
Contributions to this bulletin would not require any transfer of copyright,
and material presented here can be published elsewhere.\\[0.1cm]
\textbf{Subscription.}
To receive this bulletin (free) to your
e-mailbox, e-mail us:\\
{boaz.tsaban@weizmann.ac.il}
}

\newcommand{\nArxPaper}[5]{\subsection{#2}{#4}\par\hfill{\arx{#1}}\par\hfill\emph{#3}}
\newcommand{\DOIpaper}[5]{\subsection{#2}{#4}\par\hfill{\texttt{http://dx.doi.org/#1}}\par\hfill\emph{#3}}

\newcommand{\nAMSPaper}[4]{\subsection{#2}{#4}\par\hfill{\texttt{#1}}\par\hfill\emph{#3}}
\newcommand{\naAMSPaper}[3]{\subsection{#1}~\par\hfill{\texttt{#3}}\par\hfill\emph{#2}}

\newcommand{\arx}[1]{\texttt{http://arxiv.org/math/#1}}
\newcommand{\url}[1]{\bq\texttt{#1}\eq}
\newcommand{\online}[1]{The paper is available online at \url{#1}}

\title[$\mathcal{SPM}$ Bulletin \textbf{\issuenumber} (\issuemonth{} \issueyear)]{%
$\mathcal{SPM}$ Bulletin\\[0.5cm]
Issue number \issuenumber: \issuemonth{} \issueyear{} \CE{}}

\begin{document}
\maketitle

\tableofcontents

\section{Editor's notes}

1. At least two fundamental problems in the field of SPM are finally
solved. These and other exciting results are announced in the present issue.

2. The \emph{Mathematics Subject Classification} list is currently
being revised. Visit \url{http://www.msc2010.org/}
and make your suggestions. The web page promises that ``Mathematical Reviews and
Zentralblatt f\"ur Mathematik will carefully consider all feedback and use it in
preparing their joint MSC revision to be used starting in 2010.''

3. Recently, the online versions of papers
published by AMS journals do not contain the tex source.
Thus, henceforth I will often not give abstracts of their papers which I
announce here (unless the abstract is plain text). You will only
get the title and link to the full paper. This is not very
informative, so I urge authors to submit relevant abstracts
directly to me.

\medskip

Contributions to the next issue are, as always, welcome.

\medskip

\by{Boaz Tsaban}{boaz.tsaban@weizmann.ac.il}

\hfill \texttt{http://www.cs.biu.ac.il/\~{}tsaban}

\section{Research announcements}

\nArxPaper{math.GN/0609601}
{Selection Principles and special sets of reals: Open problems}
{Boaz Tsaban}
{This is the chapter on SPM and special sets of reals, to appear in:
\textbf{Open Problems in Topology II} (E.\ Pearl, ed.),
Elsevier Science, 2007.

We give a selection of major open
problems involving selection principles, diagonalizations, and
covering properties for sets of real numbers.}

\nArxPaper{math.LO/0609655}
{Winning the pressing down game but not Banach Mazur}
{Jakob Kellner, Matti Pauna, and Saharon Shelah}
{Let $S$ be the set of those $\alpha\in\omega_2$ that have
cofinality $\omega_1$. It is consistent relative to a measurable
that player II (the nonempty player) wins the pressing down game
of length $\omega_1$, but not the Banach Mazur game of length
$\omega+1$ (both starting with $S$).}

\subsection{Ramsey classes of topological and metric spaces}
This paper is a follow up of the author's programme of characterizing
Ramsey classes of structures by a combination of model theory
and combinatorics. This relates the classification programme for countable
homogeneous structures (of Lachlan and Cherlin) to the proof techniques
of the structural Ramsey theory. Here we consider the classes of topological
and metric spaces which recently were studied in the context of extremally
amenable groups and of the Urysohn space. We show that Ramsey classes are
essentially classes of finite objects only. While for Ramsey classes of topological
spaces we achieve a full characterization, for metric spaces this seems to be
at present an intractable problem.

\hfill\texttt{http://dx.doi.org/10.1016/j.apal.2005.07.004}

\nby{Jaroslav Ne\v{s}et\v{r}il}

\naAMSPaper{More on partitioning triples of countable ordinals}
{Albin L. Jones}
{http://www.ams.org/journal-getitem?pii=S0002-9939-06-08538-8}

\nAMSPaper{http://www.ams.org/journal-getitem?pii=S0002-9939-06-08572-8}
{Countable compact Hausdorff spaces need not be metrizable in ZF}
{Kyriakos Keremedis and Eleftherios Tachtsis}
{We show that the existence of a countable, first countable, zero-dimensional,
compact Hausdorff space which is not second countable, hence not metrizable,
is consistent with ZF. }

\nArxPaper{math.GN/0610226}
{Every topological group is a group retract of a minimal group}
{Michael Megrelishvili}
{We show that every Hausdorff topological group is a group retract of a
minimal topological group. This first was conjectured by Pestov in 1983. Our
main result leads to a solution of some problems of Arhangelskii. One of them
is the problem about representability of a group as a quotient of a minimal
group (Problem 519 in the first edition of \emph{Open Problems in Topology}). Our
approach is based on generalized Heisenberg groups and on groups arising from
group representations on Banach spaces and in bilinear mappings.
}

\nArxPaper{math.FA/0610289}
{The complexity of classifying separable Banach spaces up to isomorphism}
{Valentin Ferenczi, Alain Louveau, and Christian Rosendal}
{It is proved that the relation of isomorphism between separable Banach spaces
is a complete analytic equivalence relation, i.e., that any analytic
equivalence relation Borel reduces to it. Thus, separable Banach spaces up to
isomorphism provide complete invariants for a great number of mathematical
structures up to their corresponding notion of isomorphism. The same is shown
to hold for
\be
\itm complete separable metric spaces up to uniform homeomorphism,
\itm separable Banach spaces up to Lipschitz isomorphism, and
\itm up to (complemented) biembeddability,
\itm Polish groups up to topological isomorphism, and
\itm Schauder bases up to permutative equivalence.
\ee
Some of the
constructions rely on methods recently developed by S. Argyros and P. Dodos.
}

\newcommand{\KlO}{\ensuremath{\mathcal{O}}}
\nArxPaper{math.LO/0610306}
{Reals $n$-generic relative to some perfect tree}
{Bernard A. Anderson}
{We say that a real $X$ is $n$-generic
relative to a perfect tree $T$ if $X$ is a path through $T$ and
for all $\Sigma^0_n (T)$ sets $S$, there exists a number $k$ such
that either $X|k \in S$ or for all $\sigma \in T$ extending $X|k$
we have $\sigma \notin S$\@.  A real $X$ is $n$-generic relative
to some perfect tree if there exists such a $T$.  We first show
that for every number $n$ all but countably many reals are
$n$-generic relative to some perfect tree.  Second, we show that
proving this statement requires ZFC$^-$ + ``$\exists$ infinitely
many iterates of the power set of $\omega$''. Third, we prove that
every finite iterate of the hyperjump, $\KlO^{(n)}$, is not
2-generic relative to any perfect tree and for every ordinal
$\alpha$ below the least $\lambda$ such that $\sup_{\beta <
\lambda} (\beta$th admissible$) = \lambda$, the iterated hyperjump
$\KlO^{(\alpha)}$ is not 5-generic relative to any perfect tree.
Finally, we demonstrate some necessary conditions for reals to be
1-generic relative to some perfect tree.}

\nArxPaper{math.GN/0610432}
{Nagata's conjecture and countably compactifications in generic extensions}
{Lajos Soukup}
{Nagata conjectured that every $M$-space is homeomorphic to a closed subspace
of the product of a countably compact space and a metric space. This conjecture
was refuted by Burke and van Douwen, and A. Kato, independently.

However, we can show that there is a c.c.c.\ poset $P$ of size $2^{\omega}$
such that in $V^P$ Nagata's conjecture holds for each first countable regular
space from the ground model (i.e. if a first countable regular space $X\in V$
is an $M$-space in $V^P$ then it is homeomorphic to a closed subspace of the
product of a countably compact space and a metric space in $V^P$). In fact, we
show that every first countable regular space from the ground model has a first
countable countably compact extension in $V^P$, and then apply some results of
Morita. As a corollary, we obtain that every first countable regular space from
the ground model has a maximal first countable extension in model $V^P$.
}

\nArxPaper{math.GR/0610430}
{A Class of Groups in Which All Unconditionally Closed Sets are Algebraic}
{Ol'ga V. Sipacheva}
{It is proved that, in any subgroup of a direct product of countable groups,
the property of being an unconditionally closed set in the sense of Markov
coincides with that of being an algebraic set.}

\nArxPaper{math.FA/0610562}
{A $c_0$-saturated Banach space with no long unconditional basic sequences}
{Jordi Lopez Abad and Stevo Todorcevic}
{We present a Banach space $\mathfrak X$ with a Schauder basis of length
$\omega_1$ which is saturated by copies of $c_0$ and such that for every
closed decomposition of a closed subspace $X=X_0\oplus X_1$, either $X_0$ or
$X_1$ has to be separable. This can be considered as the non-separable
counterpart of the notion of hereditarily indecomposable space. Indeed, the
subspaces of $\mathfrak X$ have ``few operators'' in the sense that every
bounded operator $T:X \to \mathfrak{X}$ from a subspace $X$ of $\mathfrak{X}$
into $\mathfrak{X}$ is the sum of a multiple of the inclusion and a
$\omega_1$-singular operator, i.e., an operator $S$ which is not an
isomorphism on any non-separable subspace of $X$. We also show that while
$\mathfrak{X}$ is not distortable (being $c_0$-saturated), it is arbitrarily
$\omega_1$-distortable in the sense that for every $\lambda>1$ there is an
equivalent norm $\|\cdot \|$ on $\mathfrak{X}$ such that for every
non-separable subspace $X$ of $\mathfrak{X}$ there are $x,y\in S_X$ such that
$\|\cdot \| / \|\cdot \|\ge \lambda$.
}

\nArxPaper{math.FA/0610795}
{Spaces of continuous functions over Dugundji compacta}
{Taras Banakh and Wieslaw Kubis}
{We show that for every Dugundji compact $K$ the Banach space $C(K)$ is
1-Plichko and the space $P(K)$ of probability measures on $K$ is
Valdivia compact. Combining this result with the existence of a
non-Valdivia compact group, we answer a question of Kalenda.}

\nArxPaper{math.LO/0610988}
{Varia: Ideals and Equivalence Relations, beta-version}
{Vladimir Kanovei}
{We present a selection of basic results on Borel reducibility of Borel ideals and equivalence relations,
 especially those with comparably short proofs. The focal point are reducibility/irreducibility results related
 to some special equivalences like $E_0, E_1, E_2, E_3, E_\infty, Z_0,$ and Banach-induced equivalences
 $l_p,$ in particular several dichotomy theorems. The bulk of results included in the book were obtained in
 the 1990s, but some rather recent theorems are presented as well, like Rosendal's proof that Borel ideals
are cofinal within Borel equivalences of general form.}

\nArxPaper{math.GN/0611239}
{Equivariant embedding of metrizable $G$-spaces in linear $G$-spaces}
{Aasa Feragen}
{Given a Lie group $G$ we study the class $\M$ of proper metrizable
$G$-spaces with metrizable orbit spaces, and show that any
$G$-space $X \in \M$ admits a closed $G$-embedding into a convex
$G$-subset $C$ of some locally convex linear $G$-space, such that
$X$ has some $G$-neighborhood in $C$ which belongs to the class
$\M$. As corollaries we see that any $G$-ANE for $\M$ has the
$G$-homotopy type of some $G$-CW complex and that any $G$-ANR for
$\M$ is a $G$-ANE for $\M$.}

\nArxPaper{math.GN/0611353}
{Squares of Menger-bounded groups}
{Micha\l{} Machura, Saharon Shelah, and Boaz Tsaban}
{Using a portion of the Continuum Hypothesis, we prove that there
is a Menger-bounded (also called $o$-bounded) subgroup of the Baire
group $\Z^\N$, whose square is not Menger-bounded. This settles a
major open problem concerning boundedness notions for groups, and
implies that Menger-bounded groups need not be Scheepers-bounded.
}

\nAMSPaper{http://dx.doi.org/10.1016/j.topol.2005.10.014}
{$\kappa$-Fr\'echet-Urysohn property of $C_k(X)$}
{Masami Sakai}
{For a Tychonoff space $X$, we denote by $C_k(X)$ the space of all
real-valued continuous functions on $X$ with the compact open
topology. A space $X$ is said to be $\kappa$-Fr\'echet-Urysohn
if for every open subset $U$ of $X$ and every $x\in\cl{U}$,
there exists a sequence $\seq{x_n}$ in $U$
converging to $x$.
In this paper, we show that $C_k(X)$ is
$\kappa$-Fr\'echet-Urysohn iff every moving off family of compact subsets of $X$
has a countable subfamily which is strongly compact-finite.
In particular, we obtain that every stratifiable Baire space $C_k(X)$ is
an $M_1$-space.}

\nArxPaper{math.LO/0611744}
{How to drive our families mad}
{Saka\'e Fuchino, Stefan Geschke, and Lajos Soukup}
{
Given a family $F$ of pairwise almost disjoint sets on a countable
set $S$, we study families $F^+$ of maximal almost disjoint (mad)
sets extending $F$.

We define $a^+(F)$ to be the minimal
possible cardinality of $F^+\setminus F$ for such $F^+$, and
$a^+(\kappa)=\sup\{a^+(F): |F| \leq \kappa \}$. We show that all
infinite cardinal less than or equal to the continuum continuum
can be represented as $a^+(F)$ for some almost disjoint $F$ and
that the inequalities $\aleph_1=a<a^+(\aleph_1)=\fc$ and
$a=a^+(\aleph_1)<\fc$ are both consistent.

We also give a
several constructions of mad families with some additional
properties.
}

\nArxPaper{math.GN/0612148\label{HurTI}}
{Hurewicz sets of reals without perfect subsets}
{Dusan Repovs, Boaz Tsaban, and Lyubomyr Zdomskyy}
{We show that even for subsets $X$ of the real line which do not contain
perfect sets, the Hurewicz property does not imply the property $\sone(\Gamma,\Gamma)$,
asserting that for each countable family of open $\gamma$-covers of $X$, there is a choice
function whose image is a $\gamma$-cover of $X$.
This settles a problem of Just, Miller, Scheepers, and Szeptycki. Our main result
also answers a question of Bartoszy\'nski and the second author,
and implies that for $C_p(X)$, the conjunction of Sakai's strong countable fan tightness
and the Reznichenko property does not imply \Arh{}'s property $\alpha_2$.
}

\nArxPaper{math.LO/0612240}
{The spectrum of characters of ultrafilters on $\omega$}
{Saharon Shelah}
{We show the consistency of the set of regular cardinals which are the
character of some ultrafilter on omega is not convex. We also deal with the set
of $\pi\chi$-characters of ultrafilters on $\omega$.}

\nArxPaper{math.FA/0612307}
{Spaces of functions with countably many discontinuities}
{R. Haydon, A. Molto, and J. Orihuela}
{Let $\Gamma$ be a Polish space and let $K$ be a separable and pointwise
compact set of real-valued functions on $\Gamma$. It is shown that if each
function in $K$ has only countably many discontinuities then $C(K)$ may be
equipped with a $T_p$-lower semicontinuous and locally uniformly convex norm,
equivalent to the supremum norm.}

\nArxPaper{math.LO/0612353}
{Can groupwise density be much bigger than the non-dominating number?}
{Saharon Shelah}
{We prove that $\fg$ (the groupwise density number) is smaller or equal to $\fb^+$
(the successor of the minimal cardinality of a non-dominated subset of
$\NN$).}

\DOIpaper{10.1016/j.topol.2006.09.012}
{Productive local properties of function spaces}
{Francis Jordan}
{We characterize the spaces $X$ for which the space $C_p(X)$ of real
valued continuous functions with the topology of pointwise
convergence has local properties related to the preservation of
countable tightness or the Fr\'echet property in products. In
particular, we use the methods developed to construct an
uncountable subset $W$ of the real line such that the product of
$C_p(W)$ with any strongly Fr\'echet space is Fr\'echet. The example
resolves an open question.
}

\subsection{Pinning quasi orders with their endomorphisms}
Some general properties of abstract relations are closely examined. These
include generalizations of linearity, and properties based on `pinning' an
inequality by a pair of families of endomorphisms.To each property we try to
associate a canonical definition of an augmentation (or diminishment) that
augments (or diminishes) any given relation to one satisfying the desired
property. The motivation behind this study was to identify properties
distinguishing between the product ordering and the eventual dominance ordering
of the irrationals (the family of functions from $\N$ into $\N$), and furthermore to
identify their relationship as a member of a natural class of augmentations.
\url{http://homepage.univie.ac.at/James.Hirschorn/\\research/pinning/pinning.html}
\nby{James Hirschorn}

\nArxPaper{math.LO/0612636}
{A game on the universe of sets}
{Denis I. Saveliev}
{In set theory without the axiom of regularity, we consider a game in which
two players choose in turn an element of a given set, an element
of this element, etc.; a player wins if its adversary cannot make
any next move. Sets that are winning, i.e.\ have a winning strategy
for a player, form a natural hierarchy with levels indexed by
ordinals. We show that the class of hereditarily winning sets is
an inner model containing all well-founded sets, and that all four
possible relationships between the universe, the class of
hereditarily winning sets, and the class of well-founded sets are
consistent. We describe classes of ordinals for which it is
consistent that winning sets without minimal elements are exactly
in the levels indexed by ordinals of this class. For consistency
results, we propose a new method for getting non-well-founded
models. Finally, we establish a probability result by showing that
on hereditarily finite well-founded sets the first player wins
almost always.}

\nArxPaper{math.FA/0612598}
{Algebraic characterizations of measure algebras}
{Thomas Jech}
{We present necessary and sufficient conditions for the existence of a
countably additive measure on a complete Boolean algebra.}

\section{Problem of the Issue}

We recall from \cite{Pyt84} that a topological space $X$ has the
Pytkeev property, if for every $x\in X$ and a subset $A\subset X$
such that $x\in\bar{A}\setminus A$, there exists a countable
family $\mathcal N$ of infinite subsets of $A$ which forms a
$\pi$-network at $x$.

Tychonoff spaces $X$ such that $C_p(X)$ has the Pytkeev property were characterized in
\cite{Sa} and \cite{TS} as follows:
$C_p(X)$ has the Pytkeev property iff for every $\w$-shrinkable (clopen) $\w$-cover $\cU$ of
$X$, there exists a sequence $(\cU_n)_{n\in\N}$ of infinite
subsets of $\cU$ such that $\{\bigcap\cU_n:n\in\N\}$ is an
$\w$-cover of $X$.

According to a recent result of A.~Miller
\cite[Theorem~18]{TS}, the mentioned covering property  of $X$ implies
that $X$ has a strong measure zero with respect ro any
\emph{totally-bounded} metric  on $X$. So it is natural to ask
whether this covering property implies that $X$ has a strong
measure zero with respect to \emph{any} metric on it, which in the realm
of zero-dimensional spaces\footnote{The Pytkeev property of
$C_p(X)$ implies  that $X$ is zero-dimensional \cite{Sa}.}
 is known to be equivalent to the Rothberger
covering property $\sone(\cO,\cO)$ \cite{FM}. But we do not even know the answer to
the following basic question.

\begin{prob} \label{Pyt-Sch}
Does the Pytkeev property of $C_p(X)$ imply the Menger property $\sfin(\cO,\allowbreak\cO)$ of
$X$?
\end{prob}

This question is also motivated by the subsequent theorem.

\begin{thm}[Tsaban-Zdomskyy \cite{TZ}]\label{Pyt+Sch}
Let $X$ be a Tychonoff space. If $C_p(X)$ has the Pytkeev property
and $X$ satisfies $\ufin(\cO,\Omega)$, then $X$
satisfies $\ufin(\cO,\Gamma)$ as well as $\sone(\cO,\cO)$.
\end{thm}

Recall from \cite{NSW} that a zero-dimensional space $X$ has the Gerlits-Nagy
property $(\ast)$  if and only if it has the properties $\ufin(\cO,\Gamma)$ and
$\sone(\cO,\cO)$.

\nby{Lyubomyr Zdomskyy}

\ed